\documentclass[11pt]{article}
\usepackage{latexsym}
\usepackage{amssymb,amsmath,amsfonts}
\usepackage[mathscr]{eucal}
\usepackage{amscd}
\usepackage{amsfonts,cite}
\usepackage{epsfig, amsmath, amsthm, amssymb}
\usepackage{latexsym}
\usepackage{graphics,epstopdf}
\begin{document}
\begin{center}
{\Large{\bf Finding mixed families of special polynomials\\
 associated with Appell sequences $^{\star}$}}\\~~

{\bf Subuhi Khan${}^{1}${\footnote{$^{1}$Corresponding author; E-mail:~subuhi2006@gmail.com (Subuhi Khan)}}, Nusrat Raza${}^{2}${\footnote{$^{2}$E-mail:~nraza.maths@gmail.com (Nusrat Raza)}},  and Mahvish Ali${}^{3}$ {\footnote{$^{3}$This work has been done under Junior Research Fellowship (Award letter No. F1-17.1/2014-15/ MANF-2014-15-MUS-UTT-34170/(SA-III/Website)) awarded to the third author by the University Grants Commission, Government of India, New Delhi; E-mail:~mahvishali37@gmail.com (Mahvish Ali)}}}\\
\vspace{.25cm}
${}^{1,3}$Department of Mathematics, Aligarh Muslim University, Aligarh, India \\
${}^{2}$Women's College, Aligarh Muslim University, Aligarh, India \\
\end{center}

\parindent=8mm
\noindent
{\bf Abstract:}~~In this paper, certain mixed special polynomial families associated with Appell sequences are introduced and their properties are established. Further, operational rules providing connections between these families and the known special polynomials are established, which are used to derive the identities and results for the members of these new families. Determinantal definitions of the polynomials associated with Appell family is also derived. The approach presented is general.\\



\noindent
{\bf{1.~~Introduction and preliminaries}}\\

The concepts of the monomiality principle and operational techniques are used to combine the special polynomials to find the mixed special polynomials. These polynomials have many applications in different branches of mathematics. Recently, the Laguerre-Gould Hopper polynomials (LGHP) ${}_L{H}_n^{(m,r)}(x,y,z)$ are introduced in \cite{2SubAli}  which are defined by means of the generating function
$$C_0(-xt^m)\exp(yt+zt^r)=\sum_{n=0}^{\infty}{}_LH_n^{(m,r)}(x,y,z)\frac{t^n}{n!},\eqno(1.1)$$
where $C_0(x)$ denotes the Bessel-Tricomi function of order zero. The $n^{th}$-order Bessel-Tricomi functions $C_n(x)$ are specified by means of the generating function
$$\exp\left(t-\frac{x}{t}\right)=\sum_{n=0}^{\infty} C_n(x)t^n,\eqno(1.2)$$
for $t\neq 0$ and for all finite $x$ and are defined by the following series \cite [p.150]{2Dat}:
$$C_n(x)~=~x^{-\frac{n}{2}}J_n(2 \sqrt {x})~=~\sum_{k=0}^\infty \frac{(-1)^k\ x^k}{k!\ (n+k)!},~~~~n=0,1,2,\ldots\ ,\eqno(1.3)$$
with $J_n(x)$ being the ordinary cylindrical Bessel function of first kind \cite{2And}. The $0^{th}$-order Bessel-Tricomi function $C_{0}(x)$ is also given by the following operational definition:
$$C_{0}(\alpha x)=\exp \left(-\alpha D_{x}^{-1}\right)\{1\},\eqno(1.4)$$
where $D_{x}^{-1}$ denotes the inverse derivative operator and
$$D_{x}^{-n}\{1\}=\frac{x^{n}}{n!}.$$

The series definition for the LGHP ${}_L{H}_n^{(m,r)}(x,y,z)$ is given as \cite{2SubAli}:
$${}_LH_n^{(m,r)}(x,y,z)=n!\sum_{k=0}^{[\frac{n}{r}]}\frac{z^k{}_mL_{n-rk}(x,y)}{k!(n-rk)!},\eqno(1.5)$$
where ${}_mL_n(x,y)$ denotes the 2-variable generalized Laguerre polynomials (2VGLP), which are defined by the following series \cite [p.213~(27)]{2SGDat}:
$${}_mL_n(x,y)=n!\sum_{r=0}^{[\frac{n}{m}]}\frac{x^ry^{n-mr}}{(r!)^2(n-mr)!}.\eqno(1.6)$$

In view of definition (1.6), the LGHP are defined as \cite [p. 9933(2.7)]{2SubAli}:
$${}_LH_n^{(m,r)}(x,y,z)=n!\sum_{k,l=0}^{rk+ml\leq n}\frac{z^kx^ly^{n-rk-ml}}{k!(l!)^2(n-rk-ml)!}.\eqno(1.7)$$

The LGHP ${}_L{H}_n^{(m,r)}(x,y,z)$ are also defined as \cite{2SubAli}:
$${}_LH_n^{(m,r)}(x,y,z)=n!\sum_{k=0}^{[\frac{n}{m}]}\frac{x^kH_{n-mk}^{(r)}(y,z)}{(k!)^2(n-mk)!},\eqno(1.8)$$
where $H_{n}^{(r)}(y,z)$ are the Gould Hopper polynomials (GHP) \cite{2HWG} defined by
$$H_n^{(r)}(y,z)=n!\sum_{k=0}^{[\frac{n}{r}]}\frac{z^ky^{n-rk}}{k!(n-rk)!}.\eqno(1.9)$$

The operational correspondence between the LGHP ${}_L{H}_n^{(m,r)}(x,y,z)$ and the generalized Laguerre Polynomials ${}_mL_n(x,y)$ is \cite{2SubAli}:
$${}_LH_{n}^{(m,r)}(x,y,z)=\exp\left(z\frac{\partial^r}{\partial y^r}\right)\Big\{{}_mL_n(x,y)\Big\}\eqno(1.10)$$
and the correspondence between the LGHP ${}_L{H}_n^{(m,r)}(x,y,z)$ and the GHP $H_n^{(r)}(y,z)$ is \cite{2SubAli}:
$${}_LH_{n}^{(m,r)}(x,y,z)=\exp\left(D_x^{-1}\,\frac{\partial^m}{\partial y^m}\right)\Big\{H_{n}^{(r)}(y,z)\Big\}.\eqno(1.11)$$

We note that the LGHP ${}_LH_{n}^{(m,r)}(x,y,z)$ are also defined through the operational rule \cite{2SubAli}:
$${}_LH_{n}^{(m,r)}(x,y,z)=\exp\left(D_x^{-1}\frac{\partial^{m}}{\partial y^{m}}+z\frac{\partial^{r}}{\partial y^{r}}\right)\Big\{y^n\Big\}.\eqno(1.12)$$

The LGHP ${}_L{H}_n^{(m,r)}(x,y,z)$ are shown to be quasi-monomial \cite{2Dat, 2JFStef} under the action of the following multiplicative and derivative operators \cite{2SubAli}:
$$\hat{M}_{LH}:=y+mD_x^{-1}\frac{\partial^{m-1}}{\partial y^{m-1}}+rz\frac{\partial^{r-1}}{\partial y^{r-1}}\eqno(1.13)$$
and
$$\hat{P}_{LH}:=\frac{\partial}{\partial y}, \eqno(1.14)$$
respectively.\\

Consequently, $\hat{M}_{LH}$ and $\hat{P}_{LH}$ satisfy the following recurrences:

$$\hat{M}_{LH}\{{}_L{H}_n^{(m,r)}(x,y,z)\}={}_L{H}_{n+1}^{(m,r)}(x,y,z)\eqno(1.15)$$
and
$$\hat{P}_{LH}\{{}_L{H}_n^{(m,r)}(x,y,z)\}=n{}_L{H}_{n-1}^{(m,r)}(x,y,z),\eqno(1.16)$$
for all $n\in\mathbb{N}$.\\

The operators $\hat{M}_{LH}$ and $\hat{P}_{LH}$ also satisfy the following commutation relation:
$$[\hat{P}_{LH}, \hat{M}_{LH}]=1.\eqno(1.17)$$

In view of the monomiality principle equation
$$\hat{M}_{LH} \hat{P}_{LH}\{{}_L{H}_n^{(m,r)}(x,y,z)\}=n{}_L{H}_n^{(m,r)}(x,y,z),\eqno(1.18)$$
the differential equation satisfied by ${}_LH_n^{(m,r)}(x,y,z)$ is \cite{2SubAli}
$$\left(m\frac{\partial^{m}}{\partial y^{m}}+rz\frac{\partial^{r+1}}{\partial x\partial y^{r}}+y\frac{\partial^2}{\partial x\partial y}-n\frac{\partial}{\partial x}\right){}_LH_{n}^{(m,r)}(x,y,z)=0.\eqno(1.19)$$

Also, ${}_L{H}_n^{(m,r)}(x,y,z)$ can be explicitly constructed as:
$${}_L{H}_n^{(m,r)}(x,y,z)=\hat{M}_{LH}^{n}\{1\},~~~~~{}_L{H}_0^{(m,r)}(x,y,z)=1.\eqno(1.20)$$

Identity (1.20) implies that the exponential generating function of ${}_L{H}_n^{(m,r)}(x,y,z)$ can be given in the form:
$$\exp(t\hat{M}_{LH})\{1\}=\sum\limits_{n=0}^{\infty}{{}_L{H}_n^{(m,r)}(x,y,z)}\frac{t^n}{n!}, ~~~~|t|<\infty.\eqno(1.21)$$

For suitable values of the indices and variables the LGHP ${}_LH_n^{(m,r)}(x,y,z)$ give a number of other known special polynomials as its special cases. We mention these special cases in the following table:\\

\noindent
\textbf{Table 1. Special cases of the LGHP ${}_LH_n^{(m,r)}(x,y,z)$ }\\
\\
{\tiny{
\begin{tabular}{|l|l|l|l|l|}
\hline
\bf{S.} & \bf{Values of } & \bf{Relation Between the }& \bf{Name of }& \bf{Generating Functions of}  \\
\bf{No.}& \bf{the Indices} & \bf{LGHP~${}_LH_n^{(m,s)}(x,y,z)$ }&\bf{the Known}& \bf{the Known Polynomials} \\
& \bf{and Variables}& \bf{and its Special Case}&  \bf{Polynomials}&\\
\hline
 I.& $m=1,~r=2$; &  &  3-Variable&  \\
& $x\rightarrow -x$ & ${}_LH_{n}^{(1,2)}(-x,y,z)={}_LH_{n}(x,y,z)$ & Laguerre -&$C_0(xt)\exp\left(yt+zt^2\right)=\sum_{n=0}^{\infty}{}_LH_{n}(x,y,z)\frac{t^n}{n!}$\\
&&&Hermite \cite{2GGGA}&\\
\hline
II. & $m=1,~r=2$; & &2-Variable &  \\
&~$z=-\frac{1}{2}$,  &${}_LH_n^{(1,2)}(-x,y,-\frac{1}{2})={}_LH_{n}^{\ast}(x,y)$ & Laguerre - &$C_0(xt)\exp\left(yt-\frac{1}{2}t^2\right)=\sum_{n=0}^{\infty}{}_LH_{n}^{\ast}(x,y)\frac{t^n}{n!}$ \\
&$x\rightarrow -x$&&Hermite \cite{2GAMDat}&\\
\hline
 III.& $m=1,~r=2$;  & & Laguerre -&  \\
& ~$y=1$, $z\rightarrow y$, & ${}_LH_{n}^{(1,2)}(-x,1,y)=\varphi_n(x,y)$& Hermite&$C_0(xt)\exp\left(t+yt^2\right)=\sum_{n=0}^{\infty}\varphi_{n}(x,y)\frac{t^n}{n!}$ \\
&$x\rightarrow -x$&& type \cite{2GGADat}&\\
\hline
IV. & $x=0$ & ${}_LH_{n}^{(m,r)}(0,y,z)=H_n^{(r)}(y,z)$ &Gould-& $\exp\left(yt+zt^r
\right)=\sum_{n=0}^{\infty}H^{(r)}_{n}(y,z)\frac{t^n}{n!}$  \\
&&&Hopper \cite{2HWG}&\\
\hline
V.& $z=0$ &  & 2-Variable &  \\
& &${}_LH_{n}^{(m,r)}(x,y,0)={}_mL_{n}(x,y)$&generalized&$C_0\left(-xt^m\right)\exp\left(yt\right)=\sum_{n=0}^{\infty}{}_mL_n(x,y)\frac{t^n}{n!}$\\
&&&Laguerre \cite{2SGDat}& \\
\hline
VI.&   $r=m$; $x=0$, &  &  2-Variable & \\
 & $y\rightarrow -D_x^{-1}$,& ${}_LH_{n}^{(m,m)}(0,-D_x^{-1},y)$ &generalized&$C_0(xt)\exp\left(yt^m\right)=\sum_{n=0}^{\infty}{}_{[m]}L_n(x,y)\frac{t^n}{n!}$ \\
& $z\rightarrow y$&$={}_{[m]}L_n(x,y)$&Laguerre&\\
&&&type \cite{2GDat}&\\
\hline
VII. & $r=m-1;~x=0$, & ${}_LH_{n}^{(m,m-1)}(0,x,y)=U_n^{(m)}(x,y)$ &Generalized &$\exp\left(xt+yt^{m-1}\right)=\sum_{n=0}^{\infty}U_{n}^{(m)}(x,y)\frac{t^n}{n!}$  \\
&$y\rightarrow x$, $z\rightarrow y$&&Chebyshev \cite{2DAT}&\\
\hline
VIII. &$m=1$; $z=0$, & ${}_LH_{n}^{(1,r)}(-x,y,0)=L_n(x,y)$ & 2-Variable& $C_0(xt)\exp(yt)=\sum_{n=0}^{\infty}L_n(x,y)\frac{t^n}{n!}$\\
&$x\rightarrow -x$&&Laguerre \cite{2GADat}&\\
\hline
IX.& $m=1$; $z=0$,  & $n!{}_LH_{n}^{(1,r)}(y,-D_x^{-1},0)$ & 2-Variable&  \\
  & $x\rightarrow y$,&$=R_{n}(x,y)$&Legendre \cite{2GPED} &  $C_0(xt)C_0(-yt)=\sum_{n=0}^{\infty} \frac{R_{n}(x,y)}{n!}\frac{t^n}{n!}$\\
 &$y\rightarrow -D_x^{-1}$&&&\\
\hline
 X.& $x=0,~y\rightarrow x,$ &  &2-Variable & \\
&$z\rightarrow y\partial_yy$ &${}_LH_{n}^{(m,r)}(0,x,y\partial_yy)=e_n^{(r)}(x,y)$&truncated of&$\frac {1}{\left(1-yt^r\right)}\exp(xt) =\sum_{n=0}^{\infty}e_{n}^{(r)}(x,y)\frac{t^n}{n!}$  \\
&&&order $s$ \cite{2GDat}&\\
&&&(see \cite{2GMDattoli})&\\
\hline
XI.& $r=2$; $x=0$ &  &  2-Variable &   \\
& &${}_LH_{n}^{(m,2)}(0,y,z)=H_n(y,z)$&Hermite-Kamp$\acute{e}$&$\exp\left(yt+zt^2\right) =\sum_{n=0}^{\infty}H_n(y,z)\frac{t^n}{n!}$ \\
&&&de F$\acute{e}$riet \cite{2PApp}&\\
\hline
XII.& $r=2$; $x=0$,  & ${}_LH_{n}^{(m,2)}(0,D_x^{-1},y)$ & Hermite& $C_0(-xt)\exp\left(yt^2\right) =\sum_{n=0}^{\infty}G_n(x,y)\frac{t^n}{n!}$ \\
  & $y\rightarrow D_x^{-1}$, $z\rightarrow y$ & $=G_n(x,y)$& type \cite{2GDattoli} & \\
\hline
XIII.& $m=2$;~$z=0$,   &  &&  \\
  & $x\rightarrow(\frac{x^2-1}{4}),$ & ${}_LH_{n}^{(2,r)}\left(\frac{x^2-1}{4},x,0\right)=P_n(x)$ &Legendre \cite{2Rain} &$C_0\left(-\frac {\left(x^{2}-1\right)}{4}t^2\right)\exp(xt) =\sum_{n=0}^{\infty}P_n(x)\frac{t^n}{n!}$  \\
&$y\rightarrow x$ &&&\\
\hline
XIV.& $x\rightarrow y\partial_yy$, & ${}_LH_{n}^{(m,r)}(y\partial_yy,x,z)$ & 3-Variable& \\
  & $y\rightarrow x$ &$ =~H_{n}^{(r,m)}(x,y,z)$ & generalized&$\exp\left(xt+yt^m+zt^r\right)
=\sum_{n=0}^{\infty}H_{n}^{(r,m)}(x,y,z)\frac{t^n}{n!}$  \\
&&&Hermite \cite{2GPEI}&\\
\hline
XV.& $m=2,~r=3;$  & ${}_LH_{n}^{(2,3)}(z\partial_zz,x,y)$ &Bell-type \cite{2GADDat}&\\
  &$x\rightarrow z\partial_zz$,&$ =~H_{n}^{(3,2)}(x,y,z)$ &&$\exp\left(xt+yt^3+zt^2\right)
=\sum_{n=0}^{\infty}H_{n}^{(3,2)}(x,y,z)\frac{t^n}{n!}$ \\
&$y\rightarrow x, z\rightarrow y$&&&\\
\hline

\end{tabular}}}\\
\\

The correspondence given in Table 1 can be used to derive the results for the polynomials related to the LGHP ${}_LH_n^{(m,r)}(x,y,z).$\\

The Appell polynomials are very often found in different applications in pure and applied mathematics. The Appell polynomials \cite{2App} may be defined by either of the following equivalent conditions: \\

$\{A_{n}(x)\}$ $(n\in \mathbb{N}_{0})$, is an Appell set ($A_{n}$ being of degree exactly $n$) if either

\begin{enumerate}
\item[{(i)}] $\frac {d}{dx}A_{n}(x)=n~A_{n-1}(x),~~~~~~~~~~n\in \mathbb{N}_{0}$, or
\item[{(ii)}] there exists an exponential generating function of the form
$$A(t)\exp\left(xt\right)=\sum _{n=0}^{\infty }A_{n} (x)\frac{t^{n} }{n!},\eqno(1.22)$$
where $A(t)$ has (at least the formal) expansion:
$$A(t)=\sum _{n=0}^{\infty }A_{n} \frac{t^{n} }{n!},~~~~~~~~~~~~~ A_{0} \ne 0 .\eqno(1.23)$$
\end{enumerate}

Roman \cite{2SRom} characterized Appell sequences in several ways. We recall the following result \cite [Theorem 2.5.3]{2SRom}, which can be viewed as an alternate definition of Appell sequences:\\

The sequence $A_{n}(x)$ is Appell for $g(t)$, if and only if
$$\frac {1}{g(t)}\exp\left(xt\right)=\sum _{n=0}^{\infty }A_{n} (x)\frac{t^{n} }{n!},\eqno(1.24)$$
where
$$g(t)=\sum_{n=0}^{\infty}g_{n}\frac{t^{n}}{n!},~~~~~~~g_{0}\neq 0.\eqno(1.25)$$

In view of equations (1.22) and (1.24), we have
$$A(t)=\frac {1}{g(t)}.\eqno(1.26)$$

In this paper, the families of Laguerre-Gould Hopper based Appell Polynomials are introduced by using the concepts and the methods associated with monomiality principle. In Section 2, we introduce the  Laguerre-Gould Hopper based Appell polynomials (LGHAP) $_{{}_LH^{(m,r)}}A_{n}(x,y,z)$ and frame these polynomials within the context of monomiality principle formalism. Some operational representations are also derived. In Section 3, results are obtained for some members of Laguerre-Gould Hopper based Appell polynomial families. The article is concluded with the determinantal definition of the Laguerre-Gould Hopper based Appell polynomials (LGHAP) $_{{}_LH^{(m,r)}}A_{n}(x,y,z)$.  \\

\noindent
{\bf{2.~~Laguerre-Gould Hopper based Appell polynomials}}\\

To introduce the Laguerre-Gould Hopper based Appell polynomials ( LGHAP) denoted by $_{{}_LH^{(m,r)}}A_{n}(x,y,z)$, we prove the following result:\\

\noindent
{\bf{Theorem~2.1.}}  {\em The Laguerre-Gould Hopper based Appell polynomials (LGHAP) $_{{}_LH^{(m,r)}}A_{n}(x,y,z)$ are defined by the generating function}
$$\frac {1}{g(t)}C_0(-xt^m)\exp(yt+zt^r)=\sum _{n=0}^{\infty }{_{{}_LH^{(m,r)}}A_{n}(x,y,z)}\frac{t^{n} }{n!},\eqno(2.1)$$
{\em or, equivalently}
$$A(t)C_0(-xt^m)\exp(yt+zt^r)=\sum _{n=0}^{\infty }{_{{}_LH^{(m,r)}}A_{n}(x,y,z)}\frac{t^{n} }{n!}.\eqno(2.2)$$\\

\noindent
{\bf{Proof.}} Replacing $x$ in the l.h.s. and r.h.s of equation (1.24) by the multiplicative operator $\hat{M}_{LH}$ of the LGHP ${}_LH_{n}^{(m,r)}(x,y,z)$, we have
$$\frac {1}{g(t)} \exp\big(\hat {M}_{LH}t\big)=\sum _{n=0}^{\infty }A_{n}\big(\hat {M}_{LH}\big) \frac{t^{n} }{n!}.\eqno(2.3)$$

Using the expression of $\hat{M}_{LH}$ given in equation (1.13) and then decoupling the exponential operator in the l.h.s. of the resultant equation by using the Crofton-type identity \cite [p. 12]{2GPDat}

$$f\left(y+m\lambda\frac{d^{m-1}}{dy^{m-1}}\right)\big\{1\big\}=\exp\left(\lambda\frac{d^m}{dy^m}\right)\Big\{f(y)\Big\},\eqno(2.4)$$
we find
$$\frac {1}{g(t)}\exp{\Big(z\frac{\partial^{r}}{\partial y^{r}}\Big)} \exp{\left(\left(y+mD_x^{-1}\frac{\partial^{m-1}}{\partial y^{m-1}}\right)t\right)}=\sum _{n=0}^{\infty }A_{n}\Big(y+mD_x^{-1}\frac{\partial^{m-1}}{\partial y^{m-1}}+rz\frac{\partial^{r-1}}{\partial y^{r-1}}\Big) \frac{t^{n} }{n!},$$
which on further use of identity (2.4) in the l.h.s. becomes
$$\frac {1}{g(t)}\exp{\Big(z\frac{\partial^{r}}{\partial y^{r}}\Big)}\exp{\Big(D_x^{-1}\frac{\partial^{m}}{\partial y^{m}}\Big)} \exp({yt})=\sum _{n=0}^{\infty }A_{n}\Big(y+mD_x^{-1}\frac{\partial^{m-1}}{\partial y^{m-1}}+rz\frac{\partial^{r-1}}{\partial y^{r-1}}\Big) \frac{t^{n} }{n!}.\eqno(2.5)$$

Now, expanding the second exponential in the l.h.s. of equation (2.5) and using definition (1.4), we find
$$\frac {1}{g(t)}C_0\left(-xt^m\right)\exp{\Big(z\frac{\partial^{r}}{\partial y^{r}}\Big)} \exp({yt})=\sum _{n=0}^{\infty }A_{n}\Big(y+mD_x^{-1}\frac{\partial^{m-1}}{\partial y^{m-1}}+rz\frac{\partial^{r-1}}{\partial y^{r-1}}\Big) \frac{t^{n} }{n!}.\eqno(2.6)$$

Again, expanding the first exponential in the l.h.s. of equation (2.6) and denoting the resultant LGHAP in the r.h.s. by $_{{}_LH^{(m,r)}}A_{n}(x,y,z)$, that is
$$_{{}_LH^{(m,r)}}A_{n}(x,y,z)=A_{n}(\hat{M}_{LH})=A_{n}\Big(y+mD_x^{-1}\frac{\partial^{m-1}}{\partial y^{m-1}}+rz\frac{\partial^{r-1}}{\partial y^{r-1}}\Big),\eqno(2.7) $$
we get assertion (2.1). Also, in view of equation (1.26), generating function (2.1) can be expressed equivalently as equation (2.2).\\

Next, to show that the LGHAP $_{{}_LH^{(m,r)}}A_{n}(x,y,z)$ satisfy the monomiality property, we prove the following result:\\

\noindent
{\bf{Theorem~2.2.}} {\em The Laguerre-Gould Hopper based Appell polynomials $_{{}_LH^{(m,r)}}A_{n}(x,y,z)$ are quasi-monomial with respect to the following multiplicative and derivative operators:}
$$\hat {M}_{LHA}=y+mD_x^{-1}\frac{\partial^{m-1}}{\partial y^{m-1}}+rz\frac{\partial^{r-1}}{\partial y^{r-1}}-\frac {g^{\prime} \left(\partial_{y}\right)} {g\left(\partial_{y}\right)},\eqno(2.8a)$$
{\em or,}
$$\hat {M}_{LHA}=y+mD_x^{-1}\frac{\partial^{m-1}}{\partial y^{m-1}}+rz\frac{\partial^{r-1}}{\partial y^{r-1}}+\frac {A^{\prime} \left(\partial_{y}\right)}{A\left(\partial_{y}\right)}\eqno(2.8b)$$
{\em and}
$$\hat {P}_{LHA}=\partial_{y},\eqno(2.9)$$
{\em respectively, where $\partial_{y}:= \frac{\partial}{\partial y}$}.\\

\noindent
{\bf{Proof.}} Consider the following identity:
$$\partial_{y}~\left\{\exp \left(yt+zt^r\right)\right\}=t~\exp \left(yt+zt^r\right).\eqno(2.10)$$

Differentiating equation (2.3) partially with respect to $t$ and in view of relation (2.7), we find
$$\left(\hat {M}_{LH}-\frac{g'(t)}{g(t)}\right)\frac {1}{g(t)} \exp{\big(\hat {M}_{LH}t\big)}=\sum _{n=0}^{\infty }{_{{}_LH^{(m,r)}}A_{n+1}}(x,y,z)\frac{t^{n} }{n!},$$
which on using equations (1.21) and (1.1) gives
$$\left(\hat {M}_{LH}-\frac{g'(t)}{g(t)}\right)\frac {1}{g(t)}C_0\left(-xt^m\right)\exp\left(yt+zt^r\right)=\sum _{n=0}^{\infty}{_{{}_LH^{(m,r)}}A_{n+1}}(x,y,z)\frac{t^{n} }{n!}.\eqno(2.11)$$

In view of relation (2.10), the above equation becomes
$$\left(\hat {M}_{LH}-\frac{g'\left(\partial_{y}\right)}{g\left(\partial_{y}\right)}\right)\left\{\frac {1}{g(t)}C_0\left(-xt^m\right)\exp\left(yt+zt^r\right)\right\}=\sum _{n=0}^{\infty }{_{{}_LH^{(m,r)}}A_{n+1}}(x,y,z)\frac{t^{n} }{n!},\eqno(2.12)$$
which on using generating function (2.1) becomes
$$\left(\hat {M}_{LH}-\frac{g'\left(\partial_{y}\right)}{g\left(\partial_{y}\right)}\right)\left \{\sum _{n=0}^{\infty }{_{{}_LH^{(m,r)}}A_{n}(x,y,z)}\frac{t^{n} }{n!}\right \}=\sum _{n=0}^{\infty }{_{{}_LH^{(m,r)}}A_{n+1}}(x,y,z)\frac{t^{n} }{n!}.\eqno(2.13)$$

Adjusting the summation in the l.h.s. of equation (2.13) and then equating the coefficients of like powers of $t$, we find
$$\left(\hat {M}_{LH}-\frac{g'\left(\partial_{y}\right)}{g\left(\partial_{y}\right)}\right) \{_{{}_LH^{(m,r)}}A_{n}(x,y,z)\}=_{{}_LH^{(m,r)}}A_{n+1}(x,y,z),\eqno(2.14)$$
which, in view of equation (1.15) shows that the corresponding multiplicative operator for $_{{}_LH^{(m,r)}}A_{n}(x,y,z)$ is given as:
$$\hat{M}_{LHA}=\left(\hat {M}_{LH}-\frac{g'\left(\partial_{y}\right)}{g\left(\partial_{y}\right)}\right).$$

Finally, using equation (1.13) in the r.h.s of above equation, we get assertion (2.8a). Also, in view of relation (1.26) assertion (2.8a) can be expressed equivalently as equation (2.8b).\\

Again, in view of identity (2.10), we have
$$\partial_{y} \left\{\frac {1}{g(t)}C_0\left(-xt^m\right)\exp\left(yt+zt^r\right)\right\}=t \frac {1}{g(t)}C_0\left(-xt^m\right)\exp\left(yt+zt^r\right),$$
which on using generating function (2.1) becomes
$$\partial_{y}\left \{\sum _{n=0}^{\infty }{_{{}_LH^{(m,r)}}A_{n}(x,y,z)}\frac{t^{n} }{n!}\right \}=\sum _{n=1}^{\infty }{_{{}_LH^{(m,r)}}A_{n-1}(x,y,z)}\frac{t^{n} }{(n-1)!}.$$

Adjusting the summation in the l.h.s. of the above equation and then equating the coefficients of like powers of $t$, we get
$$ \partial_{y} \{_{{}_LH^{(m,r)}}A_{n}(x,y,z)\}=n~_{{}_LH^{(m,r)}}A_{n-1}(x,y,z),~~~~~~~~~~~~~n \ge 1,\eqno(2.15)$$
which in view of equation (1.16) yields assertion (2.9). \\

\noindent
{\bf{Theorem~2.3.}}~{\em The Laguerre-Gould Hopper based Appell polynomials $_{{}_LH^{(m,r)}}A_{n}(x,y,z)$ satisfy the following differential equation:}
$$\left(y\partial_{y}+mD_x^{-1}\frac{\partial^{m}}{\partial y^{m}}+rz\frac{\partial^{r}}{\partial y^{r}}-\frac {g^{\prime} \left(\partial_{y}\right)} {g\left(\partial_{y}\right)}\partial_{y}-n\right){_{{}_LH^{(m,r)}}A_{n}(x,y,z)}=0,\eqno(2.16a)$$\\
\noindent
{\em or, equivalently}
$$\left(y\partial_{y}+mD_x^{-1}\frac{\partial^{m}}{\partial y^{m}}+rz\frac{\partial^{r}}{\partial y^{r}}+\frac {A^{\prime} \left(\partial_{y}\right)}{A\left(\partial_{y}\right)}\partial_{y}-n\right){_{{}_LH^{(m,r)}}A_{n}(x,y,z)}=0.\eqno(2.16b)$$\\

\noindent
{\bf{Proof.}} Using equations (2.8a) and (2.9) in the corresponding equation (1.18) for the Laguerre-Gould Hopper based Appell polynomials $_{{}_LH^{(m,r)}}A_{n}(x,y,z)$ , we get assertion (2.16a) and similarly using equations (2.8b) and (2.9) in equation (1.18), we get assertion (2.16b).\\

\noindent
{\bf{Remark~2.1.}} The Laguerre-Gould Hopper based Appell polynomials $_{{}_LH^{(m,r)}}A_{n}(x,y,z)$ are defined by the following series
$$_{{}_LH^{(m,r)}}A_{n}(x,y,z)=\sum_{k=0}^{n}{n\choose k}{}_{L}H_{n-k}^{(m,r)}(x,y,z)A_{k},\eqno(2.17)$$
where $A_{k}$ is given by equation (1.23).\\

We have mentioned special cases of the LGHP ${}_{L}H_{n}^{(m,r)}(x,y,z)$ in Table 1. Now, for the same choice of the variables and indices the  LGHAP $_{{}_LH^{(m,r)}}A_{n}(x,y,z)$ reduce to the corresponding special case.  We mention these known and new special polynomials related to the Appell sequences in the following table:\\

\noindent
{\bf Table 2. Special cases of the LGHAP $_{{}_LH^{(m,r)}}A_{n}(x,y,z)$ .}\\

\noindent
{\tiny{
\begin{tabular}{|l|l|l|l|}
\hline
&&&\\
{\bf S.No.}&{\bf Values of the indices }&{\bf Relation between the LGHAP}&{\bf Name of the} \\
&{\bf and variables}&{\bf $_{{}_LH^{(m,r)}}A_{n}(x,y,z)$ and its special}&{\bf special polynomial}\\
&&{\bf case}&\\
\hline
&&&\\
I.&$m=1$, $r=2$; &$_{{}_LH^{(1,2)}}A_{n}(-x,y,z)=_{{}_LH}A_{n}(x,y,z)$&3-variable Laguerre-Hermite based\\
&$x\rightarrow -x$&&Appell polynomials (3VLHAP)\\
\hline
&&&\\
II.&$m=1$, $r=2$; &$_{{}_LH^{(1,2)}}A_{n}(-x,y,-\frac{1}{2})=_{{}_LH^{\star}}A_{n}(x,y)$&2-variable Laguerre-Hermite based \\
&$z=-\frac{1}{2}$, $x\rightarrow -x$&&Appell polynomials (2VLHAP)\\
\hline
&&&\\
III.&$m=1$, $r=2$;  &$_{{}_LH^{(1,2)}}A_{n}(-x,1,y)={}_{\varphi}A_{n}(x,y)$&2-variable Laguerre-Hermite type based \\
&$y=1$, $z\rightarrow y$, $x\rightarrow -x$ &&Appell polynomials (LHTAP)\\
\hline
&&&\\
IV.&$x=0$ &$_{{}_LH^{(m,r)}}A_{n}(0,y,z)=_{H^{(r)}}A_{n}(y,z)$&Gould Hopper based\\
&&& Appell polynomials (GHAP)\cite{SubNus}\\
\hline
&&&\\
V.&$z=0$&$_{{}_LH^{(m,r)}}A_{n}(x,y,0)=_{{}_mL}A_{n}(x,y)$&2-variable generalized Laguerre \\
&&&based Appell polynomials (2VGLAP)\\
\hline
&&&\\
VI.&$r=m$; $x=0$, &$_{{}_LH^{(m,m)}}A_{n}(0,-D_{x}^{-1},y)=_{{}_{[m]}L}A_{n}(x,y)$&2-variable generalized Laguerre type based \\
&$y\rightarrow -D_{x}^{-1}$, $z\rightarrow y$&&Appell polynomials (2VGLTAP)\\
\hline
&&&\\
VII.&$r=m-1$; $x=0$, &$_{{}_LH^{(m,m-1)}}A_{n}(0,x,y)=_{U^{(m)}}A_{n}(x,y)$&generalized Chebyshev based \\
&$y\rightarrow x$, $z\rightarrow y$&&Appell polynomials (GCAP)\\
\hline
&&&\\
VIII.&$m=1$;  $z=0$,&$_{{}_LH^{(1,r)}}A_{n}(-x,y,0)={}_{L}A_{n}(x,y)$&2-variable Laguerre based \\
&$x\rightarrow -x$&&Appell polynomials (2VLAP)\cite{SubSaad}\\
\hline
&&&\\
IX.&$m=1$; $z=0$,&$_{{}_LH^{(1,r)}}A_{n}(y,-D_x^{-1},0)=\frac{{}_{R}A_{n}(x,y)}{n!}$&2-variable Legendre based \\
&$x\rightarrow y$, $y\rightarrow -D_x^{-1}$&&Appell polynomials (2VLeAP)\\
\hline
&&&\\
X.&$x=0$, $y\rightarrow x$, &$_{{}_LH^{(m,r)}}A_{n}(0,x,y\partial_{y}y)=_{e^{(r)}}A_{n}(x,y)$&2-variable truncated exponential based \\
&$z\rightarrow y\partial_{y}y$&&Appell polynomials (2VTEAP)\\
\hline
&&&\\
XI.&$r=2$; $x=0$&$_{{}_LH^{(m,2)}}A_{n}(0,y,z)={}_{H}A_{n}(y,z)$&2-variable Hermite Kamp${\rm\acute e}$ de F${\acute{e}}$riet based \\
&&&Appell polynomials (2VHKdFAP)\\
\hline
&&&\\
XII.&$r=2$; $x=0$, &$_{{}_LH^{(m,2)}}A_{n}(0,y,z)={}_{G}A_{n}(x,y)$&Hermite type based \\
&$y\rightarrow D_{x}^{-1}$, $z\rightarrow y$&&Appell polynomials (HTAP)\\
\hline
&&&\\
XIII.&$m=2$; $z=0$,&$_{{}_LH^{(2,r)}}A_{n}\left(\frac {x^{2}-1}{4},x,0\right)={}_{P}A_{n}(x)$&Legendre based \\
&$x\rightarrow \left(\frac {x^{2}-1}{4}\right)$, $y\rightarrow x$&&Appell polynomials (LeAP)\\
\hline
&&&\\
XIV.&$x\rightarrow y\partial_{y}y$, $y\rightarrow x$&$_{{}_LH^{(m,r)}}A_{n}(y\partial_{y}y,x,z)=_{H^{(r,m)}}A_{n}(x,y,z)$&3-variable generalized Hermite based \\
&&&Appell polynomials (3VGHAP)\\
\hline
&&&\\
XV.&$m=2$, $r=3$; &$_{{}_LH^{(2,3)}}A_{n}(z\partial_{z}z,x,y)=_{H^{(3,2)}}A_{n}(x,y,z)$&Bell type based \\
&$x\rightarrow z\partial_{z}z$, $y\rightarrow x$, $z\rightarrow y$&&Appell polynomials (BTAP)\\
\hline
\end{tabular}}}\\

\vspace{.5cm}
\noindent
{\bf{Remark~2.2.}} Keeping in view of the special cases mentioned in Table 2, we can obtain the results for the special polynomials related to the Appell sequences.\\

Next, we derive certain operational representations for the  LGHAP $_{{}_LH^{(m,r)}}A_{n}(x,y,z)$.\\

To establish the operational representation for the  LGHAP $_{{}_LH^{(m,r)}}A_{n}(x,y,z)$, we prove the following results:\\

\noindent
{\bf{Theorem~2.4.}}~{\em The following operational representation between the LGHAP $_{{}_LH^{(m,r)}}A_{n}(x,y,z)$  and the Appell polynomials $A_{n}(x)$ holds true:}
$$_{{}_LH^{(m,r)}}A_{n}(x,y,z)=\exp\left(D_{x}^{-1}\frac{\partial^{m}}{\partial y^{m}}+z\frac{\partial^{r}}{\partial y^{r}}\right)A_{n}(y).\eqno(2.18)$$

\noindent
{\bf{Proof.}}~
In view of equation (2.7), the proof is direct use of identity (2.4).\\

\noindent
{\bf{Theorem~2.5.}}~{\em The following operational representation between the LGHAP $_{{}_LH^{(m,r)}}A_{n}(x,y,z)$  and the 2VGLAP $_{{}_mL}A_{n}(x,y)$ holds true:
$$_{{}_LH^{(m,r)}}A_{n}(x,y,z)=\exp\left(z\frac{\partial^{r}}{\partial y^{r}}\right){_{{}_mL}A_{n}(x,y)}.\eqno(2.19)$$}

\noindent
{\bf{Proof.}}~From equation (2.1) (or (2.2)), we have
$$\frac{\partial^{r}}{\partial y^{r}}{_{{}_LH^{(m,r)}}A_{n}(x,y,z)}=\frac{\partial}{\partial z}{_{{}_LH^{(m,r)}}A_{n}(x,y,z)}.\eqno(2.20)$$

Since, in view of Table 1(V), we have
$$_LH^{(m,r)}_{n}(x,y,0)={}_m{L}_{n}(x,y).\eqno(2.21)$$

Consequently, from Table 2(V), we have

$$_{{}_LH^{(m,r)}}A_{n}(x,y,0)=_{{}_mL}{A_{n}(x,y)}\eqno(2.22)$$

Now, solving equation (2.20) subject to initial condition (2.22), we get assertion (2.19).\\

\noindent
{\bf{Theorem~2.6.}}~{\em The following operational representation between the LGHAP $_{{}_LH^{(m,r)}}A_{n}(x,y,z)$  and the GHAP ${}_{H^{(r)}}A_{n}(y,z)$ holds true:
$$_{{}_LH^{(m,r)}}A_{n}(x,y,z)=\exp\left(D_{x}^{-1}\frac{\partial^{m}}{\partial y^{m}}\right){}_{H^{(r)}}A_{n}(y,z).\eqno(2.23)$$}

\noindent
{\bf{Proof.}}~From equations (1.4) and (2.1) (or (2.2)), we have
$$\frac{\partial^{m}}{\partial y^{m}}~{_{{}_LH^{(m.r)}}A_{n}(x,y,z)}=\frac{\partial}{\partial D_{x}^{-1}}~{_{{}_LH^{(m,r)}}A_{n}(x,y,z)}\eqno(2.24)$$
where (\cite{2GHMDat}; p. 32 (8)):
$$\frac{\partial}{\partial D_{x}^{-1}}:=\partial_{x}x\partial_{x}.$$

Since, in view of Table 1(IV), we have
$${}_LH_{n}^{(m,r)}(0,y,z)=H_{n}^{(r)}(y,z).\eqno(2.25)$$

Consequently, from Table 2(IV), we have
$$_{{}_LH^{(m,r)}}A_{n}(0,y,z)=_{H^{(r)}}A_{n}(y,z).\eqno(2.26)$$

Solving equation (2.24) subject to initial condition (2.26), we get assertion (2.23).\\

\noindent
{\bf{3.~~Examples}}\\

The Appell polynomials have been studied because of their remarkable applications not only in different branches of mathematics but also in theoretical physics and chemistry. These polynomials contains important sequences such as Bernoulli, Euler, Genocchi,  Miller-Lee polynomials {\it etc}.\\

We establish the series definition for the Laguerre-Gould Hopper based Appell polynomials $_{{}_LH^{(m,r)}}A_{n}(x,y,z)$ .\\

Using definition (1.24) in equation (2.1), we have
$$\exp (zt^{r})C_{0}(-xt^{m})\sum\limits_{n=0}^{\infty}A_{n}(y)\frac{t^{n}}{n!}=\sum\limits_{n=0}^{\infty}{_{{}_LH^{(m,r)}}A_{n}(x,y,z)}\frac{t^{n}}{n!},\eqno(3.1)$$
which on using definition (1.3) in l.h.s. and expanding the exponential becomes
$$\sum\limits_{n,l,k=0}^{\infty}\frac{x^{k}z^{l}A_{n}(y)t^{n+rl+mk}}{(k!)^{2}l!n!}=\sum\limits_{n=0}^{\infty}{_{{}_LH^{(m,r)}}A_{n}(x,y,z)}\frac{t^n}{n!}.$$

Equating the coefficients of like powers of $t$ in the above equation we get the following series definition for the Laguerre-Gould Hopper based Appell polynomials $_{{}_LH^{(m,r)}}A_{n}(x,y,z):$
$$_{{}_LH^{(m,r)}}A_{n}(x,y,z)=n!\sum\limits_{l=0}^{[\frac{n}{r}]}\sum\limits_{k=0}^{[\frac{n}{m}]}\frac {A_{n-rl-mk}(y)z^{l}x^{k}}{l!(k!)^{2}(n-rl-mk)!}.\eqno(3.2)$$

Next, we consider some members of the Appell family and find the results for the corresponding members of the Laguerre-Gould Hopper based Appell polynomial family $_{{}_LH^{(m,r)}}A_{n}(x,y,z)$.\\

\noindent
{\bf{Example 1.}} Since, for $g(t)=\frac {(e^{t} -1)}{t}$ (or $A(t)=\frac {t}{(e^{t} -1)}$), the Appell polynomial $A_{n}(x)$ becomes the Bernoulli polynomial $B_{n}(x)$ \cite{2Rain}. Therefore, for the same choice of $g(t)$ (or $A(t)$), the LGHAP reduce to the Laguerre-Gould Hopper based Bernoulli polynomials LGHBP $_{{}_LH^{(m,r)}}B_{n}(x,y,z)$ .\\

Thus, by taking the above value of $g(t)$ (or $A(t)$) in equations (2.1) (or (2.2)), (2.8a) (or (2.8b)) and (2.9) and (3.2), we get the following generating function, multiplicative and derivative operators and the series definition respectively for the LGHBP $_{{}_LH^{(m,r)}}B_{n}(x,y,z)$:

$$\frac{t}{(e^t-1)}C_0(-xt^m)\exp(yt+zt^r)=\sum _{n=0}^{\infty }{_{{}_LH^{(m,r)}}B_{n}(x,y,z)}\frac{t^{n} }{n!},\eqno(3.3)$$

$$\hat{M}=y+mD_x^{-1}\frac{\partial^{m-1}}{\partial y^{m-1}}+rz\frac{\partial^{r-1}}{\partial y^{r-1}}+\frac{e^{\partial_{y}}(1-\partial_{y})-1}{\partial_{y}(e^{\partial_{y}}-1)}\eqno(3.4)$$
and
$$\hat{P}=\partial_{y}\eqno(3.5)$$
and
$$_{{}_LH^{(m,r)}}B_{n}(x,y,z)=n!\sum\limits_{l=0}^{[\frac{n}{r}]}\sum\limits_{k=0}^{[\frac{n}{m}]}\frac {B_{n-rl-mk}(y)z^{l}x^{k}}{l!(k!)^{2}(n-rl-mk)!}.\eqno(3.6)$$\\
Now, we draw the surface plot of the LGHBP $_{{}_LH^{(m,r)}}B_{n}(x,y,z)$. To draw the surface plot of LGHBP, we consider the values of the first six Bernoulli polynomials $B_{n}(x)$ in Table 3.\\

\noindent
{\bf Table 3. First six Bernoulli polynomials}\\

\noindent
{\tiny{
\begin{tabular}{|l|l|l|l|l|l|l|}
\hline
&&&&&&\\
$n$ & 0 & 1 & 2 & 3 & 4 & 5  \\
\hline
&&&&&&\\
$B_n(x)$ &  1  &  $x-\frac{1}{2}$  & $x^2-x+\frac{1}{6}$  &   $x^3-\frac{3}{2}x^{2}+\frac{x}{2}$  &   $x^{4}-2x^{3}+x^{2}-\frac{1}{30}$  &  $x^{5}-\frac{5}{2}x^{4}+\frac{5}{3}x^{3}-\frac{x}{6}$\\
&&&&&&\\
\hline
\end{tabular}}}\\
\\

Next, we find the value of the LGHBP $_{{}_LH^{(m,r)}}B_{n}(x,y,z)$ for $n=4$, $m=3$ and $r=5$ from equation (3.6), so that we have
$$_{{}_LH^{(3,5)}}B_{4}(x,y,z)=B_{4}(y)+24B_{1}(y)x.\eqno(3.7)$$

Using the particular values of $B_{n}(x)$ from Table 3 in equation (3.7), we find
$$_{{}_LH^{(3,5)}}B_{4}(x,y,z)=y^{4}-2y^{3}+y^{2}+24xy-12x-\frac{1}{30}.\eqno(3.8)$$

In view of equation (3.8), we  get the following surface plot of $_{{}_LH^{(3,5)}}B_{4}(x,y,z)$:\\

\begin{figure}[htb]
\begin{center}

\epsfig{file=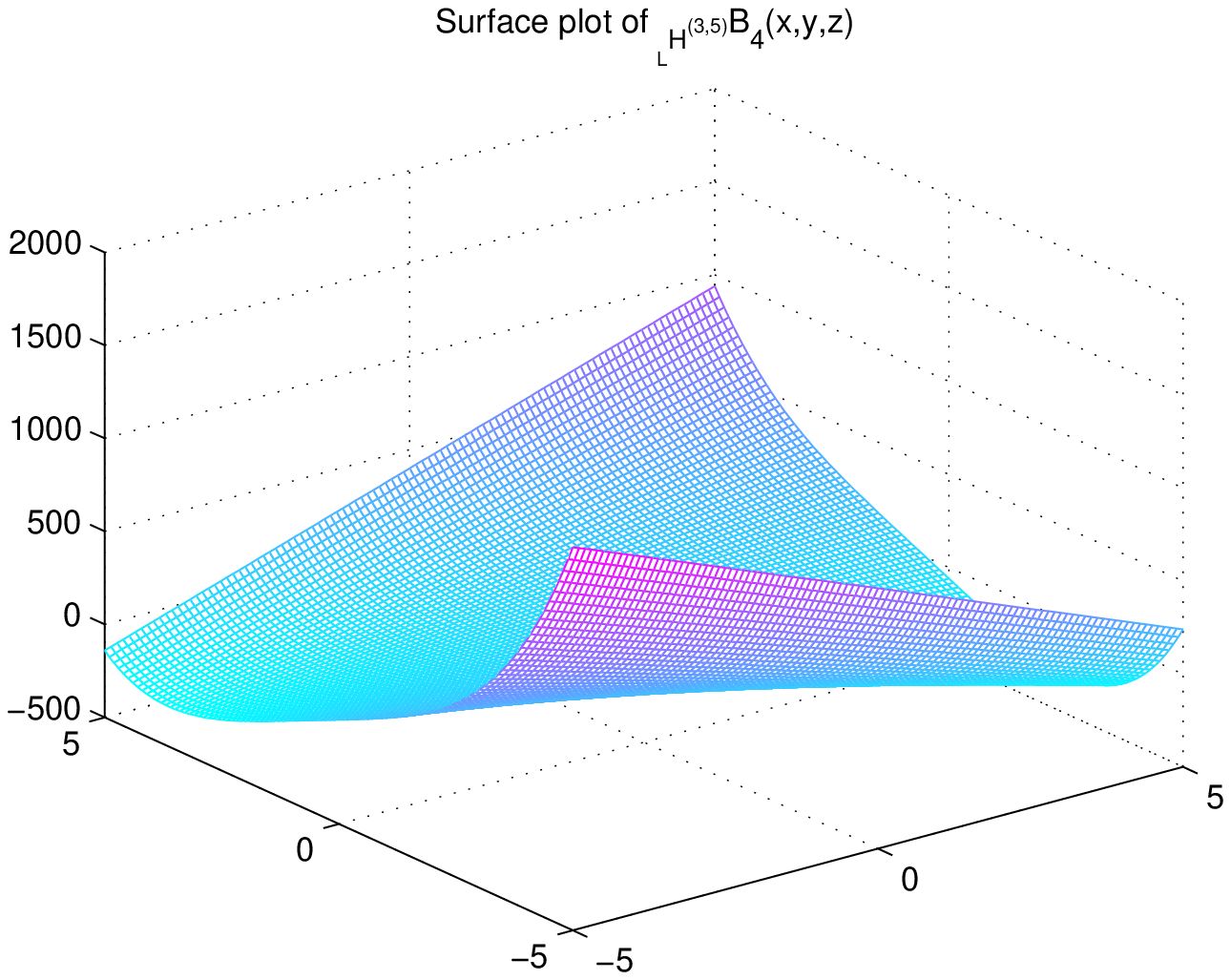, width=10cm}
\end{center}

\end{figure}

\noindent
{\bf{Example 2.}} Since, for $g(t)=\frac {(e^{t} +1)}{2}$ (or $A(t)=\frac{2}{(e^{t} +1)}$), the Appell polynomial $A_{n}(x)$ becomes the Euler polynomial $E_{n}(x)$ \cite{2Rain}. Therefore, for the same choice of $g(t)$ (or $A(t)$) the LGHAP reduce to the Laguerre-Gould Hopper based Euler polynomials LGHEP $_{{}_LH^{(m,r)}}E_{n}(x,y,z)$ .\\

Thus, by taking the above value of $g(t)$ (or $A(t)$) in equations (2.1) (or (2.2)), (2.8a) (or (2.8b)) and (2.9) and (3.2), we get the following generating function, multiplicative and derivative operators and the series definition respectively for the LGHEP $_{{}_LH^{(m,r)}}E_{n}(x,y,z)$:

$$\frac{2}{(e^t+1)}C_0(-xt^m)\exp(yt+zt^r)=\sum _{n=0}^{\infty }{_{{}_LH^{(m,r)}}E_{n}(x,y,z)}\frac{t^{n} }{n!},\eqno(3.9)$$

$$\hat{M}=y+mD_x^{-1}\frac{\partial^{m-1}}{\partial y^{m-1}}+rz\frac{\partial^{r-1}}{\partial y^{r-1}}-\frac{e^{\partial_{y}}}{(e^{\partial_{y}}+1)}\eqno(3.10)$$
and
$$\hat{P}=\partial_{y}\eqno(3.11)$$
and
$$_{{}_LH^{(m,r)}}E_{n}(x,y,z)=n!\sum\limits_{l=0}^{[\frac{n}{r}]}\sum\limits_{k=0}^{[\frac{n}{m}]}\frac {E_{n-rl-mk}(y)z^{l}x^{k}}{l!(k!)^{2}(n-rl-mk)!}.\eqno(3.12)$$

Now, we draw the surface plot of the LGHEP $_{{}_LH^{(m,r)}}E_{n}(x,y,z)$. To draw the surface of LGHEP, we consider the values of the first six Euler polynomials $E_{n}(x)$ in Table 4.\\

\noindent
{\bf{Table 4. First six Euler polynomials}}\\

\noindent
{\tiny{
\begin{tabular}{|l|l|l|l|l|l|l|}
\hline
&&&&&&\\
$n$ & 0 & 1 & 2 & 3 & 4 & 5  \\
\hline
&&&&&&\\
$E_n(x)$ &  1  &  $x-\frac{1}{2}$  & $x^2-x$  &   $x^3-\frac{3}{2}x^{2}+\frac{1}{6}$  &   $x^{4}-2x^{3}+\frac{2}{3}x$  &  $x^{5}-\frac{5}{2}x^{4}+\frac{5}{3}x^{2}-\frac{1}{2}$\\
&&&&&&\\
\hline
\end{tabular}}}\\
\\

Next, we find the value of the LGHEP $_{{}_LH^{(m,r)}}E_{n}(x,y,z)$ for $n=4$, $m=3$ and $r=5$ from equation (3.12), so that we have
$$_{{}_LH^{(3,5)}}E_{4}(x,y,z)=E_{4}(y)+24E_{1}(y)x.\eqno(3.13)$$

Using the particular values of $E_{n}(x)$ from Table 4 in equation (3.13), we find
$$_{{}_LH^{(3,5)}}E_{4}(x,y,z)=y^{4}-2y^{3}+24xy+\frac{2}{3}y-12x.\eqno(3.14)$$

In view of equation (3.14), we  get the following surface plot of $_{{}_LH^{(3,5)}}E_{4}(x,y,z)$:\\
\newpage
\begin{figure}[htb]
\begin{center}

\epsfig{file=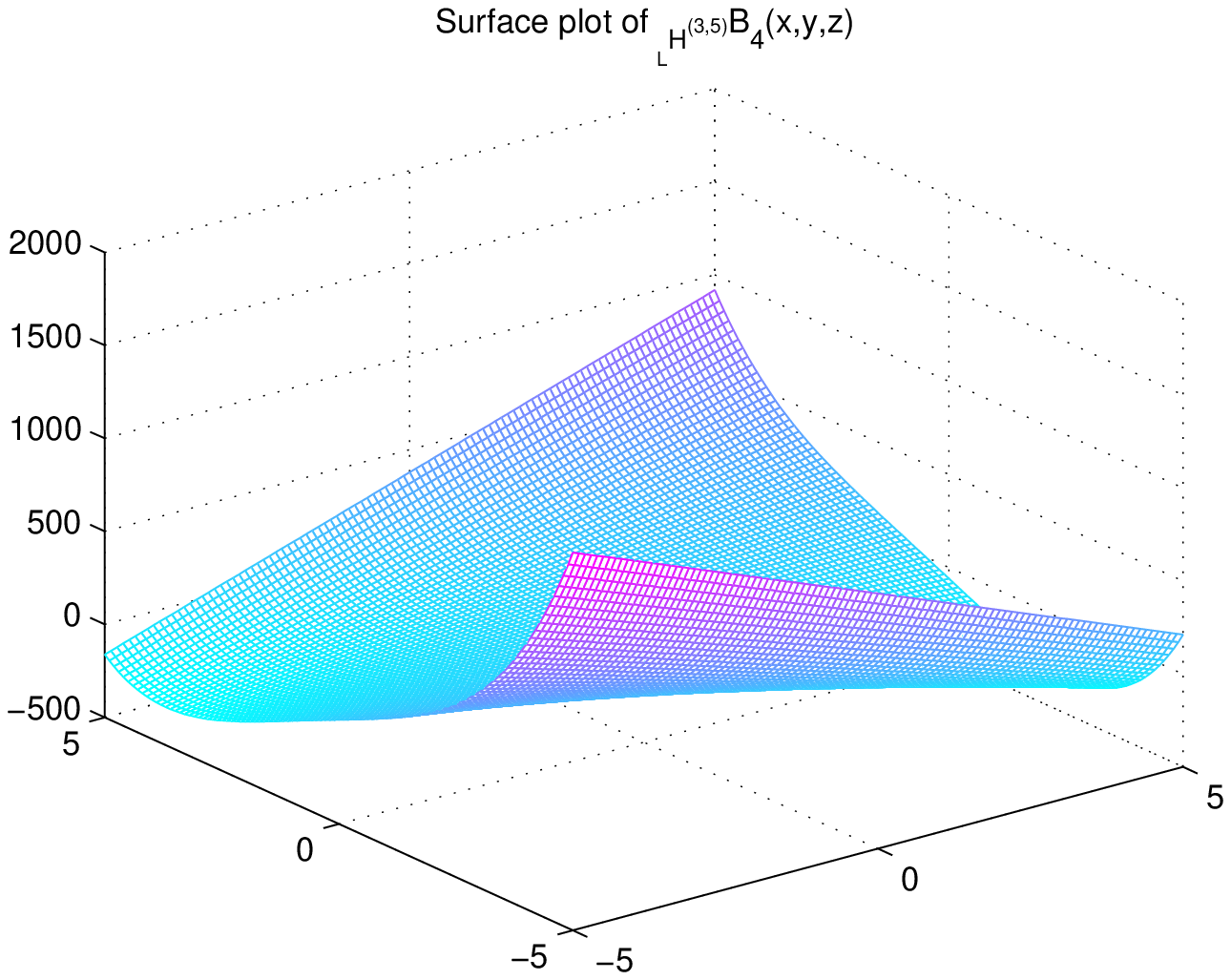, width=10cm}
\end{center}

\end{figure}

\noindent
{\bf{Example 3.}} Since, for $g(t)=(1-t)^{s+1}$  (or $A(t)=\frac{1}{(1-t)^{s+1}}$), the Appell polynomial $A_{n}(x)$ becomes the Miller-Lee polynomials $G_{n}^{(s)}(x)$ \cite{2And,Lorren} . Therefore, for the same choice of $g(t)$ (or $A(t)$) the LGHAP reduces to the Laguerre-Gould Hopper based Miller-Lee polynomials LGHMLP $_{{}_LH^{(m,r)}}G_{n}^{(s)}(x,y,z)$ .\\

Thus, by taking the above value of $g(t)$ (or $A(t)$) in equations (2.1) (or (2.2)), (2.8a) (or (2.8b)) and (2.9) and (3.2), we get the following generating function, multiplicative and derivative operators and the series definition respectively for the LGHMLP $_{{}_LH^{(m,r)}}G_{n}^{(s)}(x,y,z)$:

$$\frac{1}{(1-t)^{s+1}}C_0(-xt^m)\exp(yt+zt^r)=\sum _{n=0}^{\infty }{_{{}_LH^{(m,r)}}{G}_{n}^{(s)}(x,y,z)}t^n,\eqno(3.15)$$

$$\hat{M}=y+mD_x^{-1}\frac{\partial^{m-1}}{\partial y^{m-1}}+rz\frac{\partial^{r-1}}{\partial y^{r-1}}+(s+1)\frac{1}{1-\partial_{y}}\eqno(3.16)$$
and
$$\hat{P}=\partial_{y}\eqno(3.17)$$
and
$$_{{}_LH^{(m,r)}}G_{n}^{(s)}(x,y,z)=n!\sum\limits_{l=0}^{[\frac{n}{r}]}\sum\limits_{k=0}^{[\frac{n}{m}]}\frac {G_{n-rl-mk}^{(s)}(y)z^{l}x^{k}}{l!(k!)^{2}(n-rl-mk)!}.\eqno(3.18)$$\\

It should be noted that for $s=0$ and $s=\beta -1$, the Miller-Lee polynomials reduce to the truncated exponential polynomials $e_{n}(x)$ \cite{2And} and the modified Laguerre polynomials $f_{n}^{(\beta)}(x)$ \cite{Mag1} respectively. Therefore, for $s=0$ the LGHMLP $_{{}_LH^{(m,r)}}G_{n}^{(s)}(x,y,z)$ reduce to the Laguerre-Gould Hopper based truncated exponential polynomials LGHTP $_{{}_LH^{(m,r)}}e_{n}(x,y,z)$ and for $s=\beta-1$ the LGHMLP $_{{}_LH^{(m,r)}}G_{n}^{(s)}(x,y,z)$ reduce to the Laguerre-Gould Hopper based modified Laguerre polynomials LGHmLP $_{{}_LH^{(m,r)}}f_{n}^{(\beta)}(x,y,z)$. Thus, by taking same values of $s$ in equations (3.15)-(3.18), we can find the corresponding results for the LGHTP $_{{}_LH^{(m,r)}}e_{n}(x,y,z)$ and LGHmLP $_{{}_LH^{(m,r)}}f_{n}^{(\beta)}(x,y,z)$.\\

Now, we draw the surface plot of the LGHTP $_{{}_LH^{(m,r)}}e_{n}(x,y,z)$.\\

We have the following series definition for the LGHTP $_{{}_LH^{(m,r)}}e_{n}(x,y,z)$:
$$_{{}_LH^{(m,r)}}e_{n}(x,y,z)=n!\sum\limits_{l=0}^{[\frac{n}{r}]}\sum\limits_{k=0}^{[\frac{n}{m}]}\frac {e_{n-rl-mk}(y)z^{l}x^{k}}{l!(k!)^{2}(n-rl-mk)!}.\eqno(3.19)$$

To draw the surface plot of LGHTP, we consider the values of the first six truncated exponential polynomials $e_{n}(x)$ in Table 5.\\

\noindent
{\bf{Table 5. First six truncated exponential polynomials}}\\

\noindent
{\tiny{
\begin{tabular}{|l|l|l|l|l|l|l|}
\hline
&&&&&&\\
$n$ & 0 & 1 & 2 & 3 & 4 & 5  \\
\hline
&&&&&&\\
$e_{n}(x)$& 1 & $x+1$ & $\frac{1}{2}x^2+x+1$ & $\frac{1}{6}x^3+\frac{1}{2}x^2+x+1$ & $\frac{1}{24}x^4+\frac{1}{6}x^3+\frac{1}{2}x^2+x+1$& $\frac{1}{120}x^5+\frac{1}{24}x^4+\frac{1}{6}x^3+\frac{1}{2}x^2+x+1$\\
&&&&&&\\
\hline
\end{tabular}}}\\
\\

Next, we find the value of the LGHTP $_{{}_LH^{(m,r)}}e_{n}(x,y,z)$ for $n=4$, $m=3$ and $r=5$ from equation (3.19), so that we have
$$_{{}_LH^{(3,5)}}e_{4}(x,y,z)=e_{4}(y)+24e_{1}(y)x.\eqno(3.20)$$

Using the particular values of $e_{n}(x)$ from Table 5 in equation (3.20), we find
$$_{{}_LH^{(3,5)}}e_{4}(x,y,z)=\frac{1}{24}y^4+\frac{1}{6}y^3+\frac{1}{2}y^2+24xy+y+24x+1.\eqno(3.21)$$

In view of equation (3.21), we  get the following surface plot of $_{{}_LH^{(3,5)}}e_{4}(x,y,z)$:\\

\begin{figure}[htb]

\begin{center}
\epsfig{file=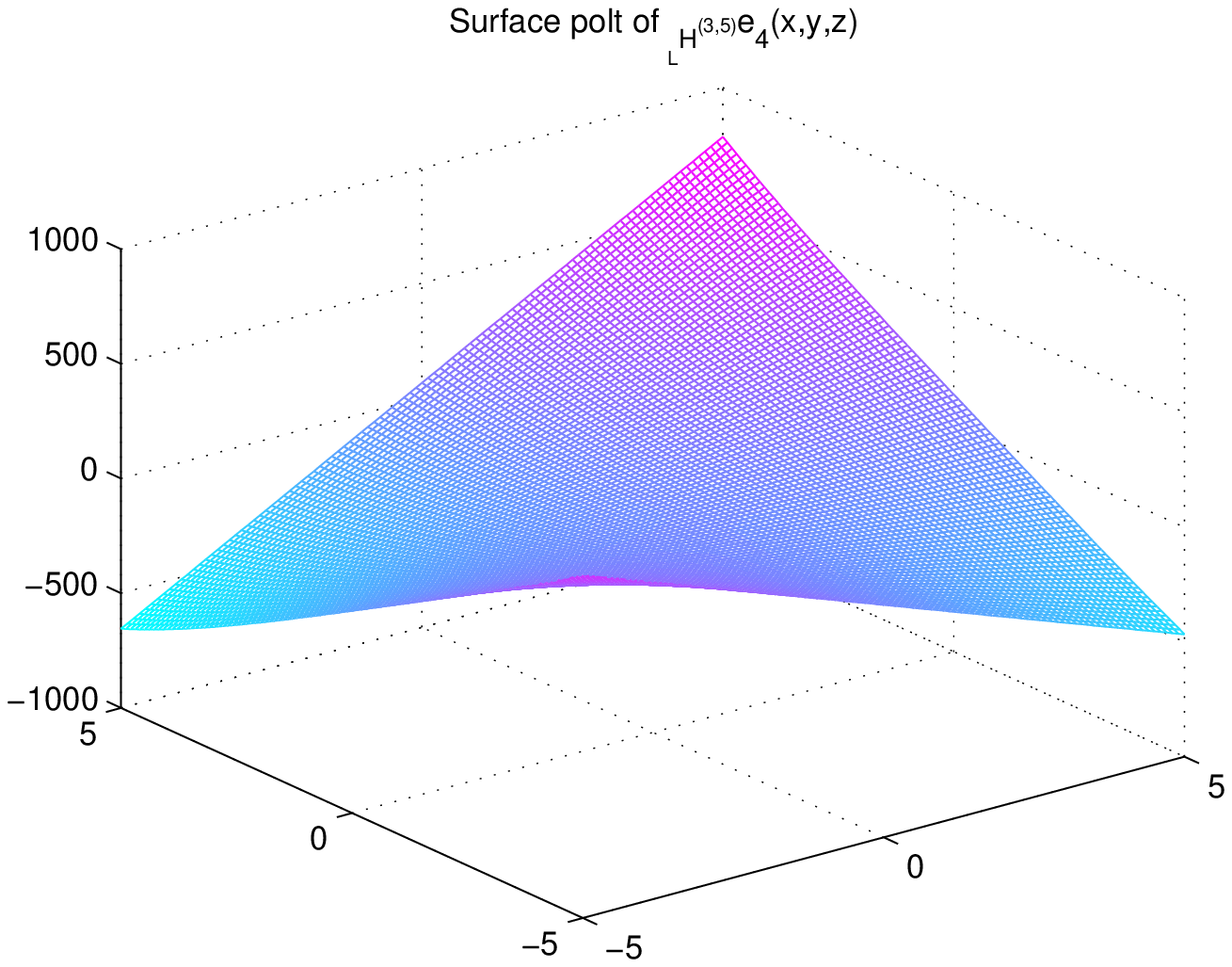, width=10cm}
\end{center}
\end{figure}

\noindent
{\bf{Example 4.}} Since, for $g(t)=\frac{(e^{t} +1)}{2t}$ (or $A(t)=\frac{2t}{(e^{t}+1)}$) the Appell polynomial $A_{n}(x)$ becomes the Genocchi polynomials  polynomial $G_{n}(x)$ \cite{2GGDat} . Therefore, for the same choice of $g(t)$ (or $A(t)$) the LGHAP reduces to the Laguerre-Gould Hopper based Genocchi polynomials LGHGP $_{{}_LH^{(m,r)}}G_{n}(x,y,z)$ .\\

Thus, by taking the above value of $g(t)$ (or $A(t)$) in equations (2.1) (or (2.2)), (2.8a) (or (2.8b)) and (2.9) and (3.2), we get the following generating function, multiplicative and derivative operators and the series definition respectively for the LGHGP $_{{}_LH^{(m,r)}}G_{n}(x,y,z)$:

$$\frac{2t}{(e^t+1)}C_0(-xt^m)\exp(yt+zt^r)=\sum _{n=0}^{\infty }{_{{}_LH^{(m,r)}}{G}_{n}(x,y,z)}\frac{t^n}{n!},\eqno(3.22)$$

$$\hat{M}=y+mD_x^{-1}\frac{\partial^{m-1}}{\partial y^{m-1}}+rz\frac{\partial^{r-1}}{\partial y^{r-1}}+\frac{e^{\partial_{y}}(1-{\partial_{y}})+1}{(e^{\partial_{y}}+1)}\eqno(3.23)$$
and
$$\hat{P}=\partial_{y}\eqno(3.24)$$
and
$$_{{}_LH^{(m,r)}}G_{n}(x,y,z)=n!\sum\limits_{l=0}^{[\frac{n}{r}]}\sum\limits_{k=0}^{[\frac{n}{m}]}\frac {G_{n-rl-mk}(y)z^{l}x^{k}}{l!(k!)^{2}(n-rl-mk)!}.\eqno(3.25)$$\\

Now, we draw the surface plot of the LGHGP $_{{}_LH^{(m,r)}}G_{n}(x,y,z)$. To draw the surface plot of LGHGP, we consider the values of the first six Genocchi polynomials $G_{n}(x)$ in Table 6.\\

\vspace{0.5cm}
\noindent
{\bf Table 6. First six Genocchi polynomials}\\

\noindent
{\tiny{
\begin{tabular}{|l|l|l|l|l|l|l|}
\hline
&&&&&&\\
$n$ & 0 & 1 & 2 & 3 & 4 & 5  \\
\hline
&&&&&&\\
$G_{n}(x)$& 0 & 1 & $2x-1$ & $3x^2-3x$ & $4x^3-6x^2+1$ & $5x^4-10x^3+5x$\\
&&&&&&\\
\hline
\end{tabular}}}\\
\\

Next, we find the value of the LGHGP $_{{}_LH^{(m,r)}}G_{n}(x,y,z)$ for $n=4$, $m=3$ and $r=5$ from equation (3.25), so that we have
$$_{{}_LH^{(3,5)}}G_{4}(x,y,z)=G_{4}(y)+24G_{1}(y)x.\eqno(3.26)$$

Using the particular values of $G_{n}(x)$ from Table 6 in equation (3.26), we find
$$_{{}_LH^{(3,5)}}G_{4}(x,y,z)=4y^3-6y^2+24x+1.\eqno(3.27)$$

In view of equation (3.27), we  get the following surface plot of $_{{}_LH^{(3,5)}}G_{4}(x,y,z)$:\\

\begin{figure}[htb]

\begin{center}
\epsfig{file=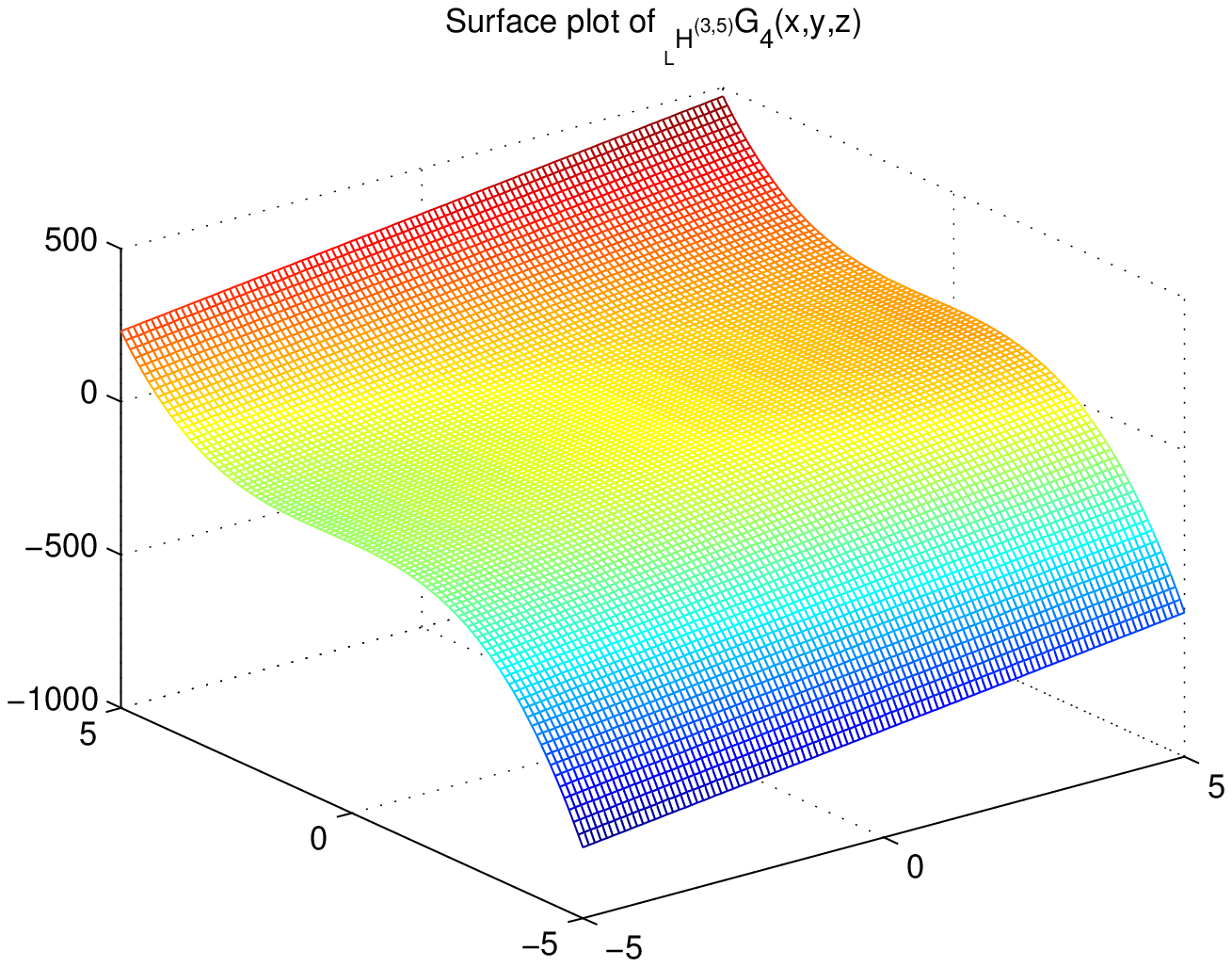, width=10cm}
\end{center}
\end{figure}

Similarly, for the other members of the Appell family, there exists a new special polynomial belonging to the LGHA family. The generating function and other properties of these special polynomials can be obtained from the results derived in Section 2.\\

We present these results along with the name and notation of the resultant special polynomial belonging to the LGHAP family in Table 7.\\

\newpage
\noindent
{\bf Table 7.  Certain results for the members belonging to the LGHAP family}\\
\\
{\tiny{
\begin{tabular}{|l|l|l|l|}
\hline
&&&\\
\bf{S. No.}  & \bf{Name/Notation of the }& \bf{Generating Function}&{\bf Multiplicative and Derivative Operators}\\
&{\bf Resultant Special}&& \\
&{\bf Polynomial}&&\\
\hline
&&&\\
I.&$_{{}_LH^{(m,r)}}B_{n}^{(\alpha)}(x,y,z)$:=&$\frac{t^{\alpha}}{(e^{t}-1)^{\alpha}}C_0(-xt^m)\exp(yt+zt^r)$&$\hat{M}=y+mD_x^{-1}\frac{\partial^{m-1}}{\partial y^{m-1}}+rz\frac{\partial^{r-1}}{\partial y^{r-1}}+\frac{\alpha (e^{\partial_{y}}(1-\partial_{y})-1)}{\partial_{y}(e^{\partial_{y}}-1)}$\\
&Laguerre-Gould Hopper based&$=\sum _{n=0}^{\infty }{_{{}_LH^{(m,r)}}B_{n}^{(\alpha)}(x,y,z)}\frac{t^{n} }{n!}$&$\hat{P}=\frac{\partial}{\partial y}$\\
&generalized Bernoulli polynomials&&\\
\hline
&&&\\
II.&$_{{}_LH^{(m,r)}}E_{n}^{(\alpha)}(x,y,z)$:=&$\frac{2^{\alpha}}{(e^t+1)^{\alpha}}C_0(-xt^m)\exp(yt+zt^r)$&$\hat{M}=y+mD_x^{-1}\frac{\partial^{m-1}}{\partial y^{m-1}}+rz\frac{\partial^{r-1}}{\partial y^{r-1}}-\frac{\alpha e^{\partial_{y}}}{(e^{\partial_{y}}+1)}$\\
&Laguerre-Gould Hopper based&$=\sum _{n=0}^{\infty }{_{{}_LH^{(m,r)}}E_{n}^{(\alpha)}(x,y,z)}\frac{t^{n} }{n!}$&$\hat{P}=\frac{\partial}{\partial y}$\\
&generalized Euler polynomials&&\\
\hline
&&&\\
III.&$_{{}_LH^{(m,r)}}\mathfrak{B}_{n}^{(\alpha)}(x,y,z;\lambda)$:=&$\left(\frac{t}{\lambda e^t-1}\right)^{\alpha}C_0(-xt^m)\exp(yt+zt^r)$&$\hat{M}=y+mD_x^{-1}\frac{\partial^{m-1}}{\partial y^{m-1}}+rz\frac{\partial^{r-1}}{\partial y^{r-1}}+\frac{\alpha (\lambda e^{\partial_{y}}(1-{\partial_{y}})-1)}{{\partial_{y}}(\lambda e^{\partial_{y}}-1)}$\\
&Laguerre-Gould Hopper based&$=\sum _{n=0}^{\infty }{_{{}_LH^{(m,r)}}\mathfrak{B}_{n}^{(\alpha)}(x,y,z;\lambda)}\frac{t^{n} }{n!}$&$\hat{P}=\frac{\partial}{\partial y}$\\
&Apostol-Bernoulli polynomials of order $\alpha$&&\\
\hline
&&&\\
IV.&$_{{}_LH^{(m,r)}}\mathfrak{B}_{n}(x,y,x;\lambda)$:=&$\frac{t}{(\lambda e^t-1)}C_0(-xt^m)\exp(yt+zt^r)$&$\hat{M}=y+mD_x^{-1}\frac{\partial^{m-1}}{\partial y^{m-1}}+rz\frac{\partial^{r-1}}{\partial y^{r-1}}+\frac{(\lambda e^{\partial_{y}}(1-{\partial_{y}})-1)}{{\partial_{y}}(\lambda e^{\partial_{y}}-1)}$\\
&Laguerre-Gould Hopper based&$=\sum _{n=0}^{\infty }{_{{}_LH^{(m,r)}}\mathfrak{B}_{n}(x,y,z;\lambda)}\frac{t^{n} }{n!}$&$\hat{P}=\frac{\partial}{\partial y}$\\
&Apostol-Bernoulli polynomials&&\\
\hline
&&&\\
V.&$_{{}_LH^{(m,r)}}\mathcal{E}_{n}^{(\alpha)}(x,y,z;\lambda)$:=&$\left(\frac{2}{\lambda e^t+1}\right)^{\alpha}C_0(-xt^m)\exp(yt+zt^r)$&$\hat{M}=y+mD_x^{-1}\frac{\partial^{m-1}}{\partial y^{m-1}}+rz\frac{\partial^{r-1}}{\partial y^{r-1}}-\frac{\alpha \lambda e^{\partial_{y}}}{(\lambda e^{\partial_{y}}+1)}$\\
&Laguerre-Gould Hopper based&$=\sum _{n=0}^{\infty }{_{{}_LH^{(m,r)}}\mathcal{E}_{n}^{(\alpha)}(x,y,z;\lambda)}\frac{t^{n} }{n!}$&$\hat{P}=\frac{\partial}{\partial y}$\\
&Apostol-Euler polynomials of order $\alpha$&&\\
\hline
&&&\\
VI.&$_{{}_LH^{(m,r)}}\mathcal{E}_{n}(x,y,z;\lambda)$:=&$\frac{2}{(\lambda e^t+1)}C_0(-xt^m)\exp(yt+zt^r)$&$\hat{M}=y+mD_x^{-1}\frac{\partial^{m-1}}{\partial y^{m-1}}+rz\frac{\partial^{r-1}}{\partial y^{r-1}}-\frac{\lambda e^{\partial_{y}}}{(\lambda e^{\partial_{y}}+1)}$\\
&Laguerre-Gould Hopper based&$=\sum _{n=0}^{\infty }{_{{}_LH^{(m,r)}}\mathcal{E}_{n}(x,y,z;\lambda)}\frac{t^{n} }{n!}$&$\hat{P}=\frac{\partial}{\partial y}$\\
&Apostol-Euler polynomials&&\\
\hline
\end{tabular}}}\\
\\

\noindent
{\bf{4.~~Concluding Remarks}}\\

F.A. Costabile and E. Longo \cite{Cos2} give a new definition for Appell polynomials by means of a determinant. According to \cite[p.1533]{Cos2}, the determinantal definition of the Appell polynomials is given as:

$$\begin{array}{l}
A_0(x)=\frac{1}{\beta}_0,~~{\beta}_0=\frac{1}{A_0},\hspace{11cm}(4.1)\\

A_n(x)=\frac{{(-1)}^n}{{(\beta_0)}^{n+1}}\left|\begin{array}{ccccccc}
 1  &  x  &  x^2  & \cdots &  x^{n-1}  & x^n \\
 \\
 {\beta}_0  &  {\beta}_1  &  {\beta}_2  & \cdots &  {\beta}_{n-1}  &  {\beta}_n \\
 \\
 0  &  {\beta}_0  &  {2  \choose  1}{\beta}_1  & \cdots & {n-1 \choose 1}{\beta}_{n-2} & {n \choose 1}{\beta}_{n-1} \\
 \\
 0 & 0 & {\beta}_0 & \cdots & {n-1 \choose 2}{\beta}_{n-3} &  {n \choose 2}{\beta}_{n-2}\\
 . & . & . & \cdots & . & . \\
 . & . & . & \cdots & . & . \\
 0 & 0 & 0 & \cdots & {\beta}_0 & {n \choose n-1}{\beta}_1
\end{array} \right|,~{\beta}_n=-\frac{1}{A_0}\Big(\sum\limits_{k=1}^{n}{n \choose k}A_k~{\beta}_{n-k}\Big),\\
\hspace{4.35in}n=1,2,3,\cdots,\end{array}\eqno(4.2)$$
where ${\beta}_0,~{\beta}_1,\cdots,{\beta}_n \in \mathbb{R}$,~~${\beta}_0 \neq 0$.\\
\vspace{.15cm}

In this Section, we give the determinantal definition of the Laguerre-Gould Hopper based Appell polynomials (LGHAP) $_{{}_LH^{(m,r)}}A_{n}(x,y,z)$.\\

On putting $n=0$ in series definition (2.17) of the Laguerre-Gould Hopper based Appell polynomials and then using equation (4.1), we obtain:
$$_{{}_LH^{(m,r)}}A_{0}(x,y,z)=\frac{1}{\beta}_0,~~{\beta}_0=\frac{1}{A_0}.\eqno(4.3)$$

Next, we expand the determinant given in equation (4.2) with respect to the first row, so that we have
$$ \begin{array}{l}
A_n(x)=\frac{(-1)^n}{({\beta}_0)^{n+1}}\left| \begin{array}{ccccc}
                                                  {\beta}_1 & {\beta}_2 & \cdots & {\beta}_{n-1} & {\beta}_n \\
                                                  \\
                                                  {\beta}_0 & {2 \choose 1}{\beta}_1 & \cdots & {n-1 \choose 1}{\beta}_{n-2}& {n \choose 1}{\beta}_{n-1} \\
                                                  \\
                                                  0 & {\beta}_0 & \cdots & {n-1 \choose 2}{\beta}_{n-3} & {n \choose 2}{\beta}_{n-2} \\
                                                  . & . & \cdots & . & . \\
                                                  . & . & \cdots & . & . \\
                                                  0 & 0 & \cdots & {\beta}_0 & {n \choose n-1}{\beta}_1
                                                  \end{array} \right|\end{array}$$\\
\\
$$\begin{array}{l}
-\frac{(-1)^n x}{({\beta}_0)^{n+1}}\left| \begin{array}{ccccc}
                                                  {\beta}_0 & {\beta}_2 & \cdots & {\beta}_{n-1} & {\beta}_n \\
                                                  \\
                                                  0 & {2 \choose 1}{\beta}_1 & \cdots & {n-1 \choose 1}{\beta}_{n-2}& {n \choose 1}{\beta}_{n-1} \\
                                                  \\
                                                  0 & {\beta}_0 & \cdots & {n-1 \choose 2}{\beta}_{n-3} & {n \choose 2}{\beta}_{n-2} \\
                                                  . & . & \cdots & . & . \\
                                                  . & . & \cdots & . & . \\
                                                  0 & 0 & \cdots & {\beta}_0 & {n \choose n-1}{\beta}_1
                                                \end{array} \right|

+\frac{(-1)^n x^2}{({\beta}_0)^{n+1}}\left| \begin{array}{ccccc}
                                                  {\beta}_0 & {\beta}_1 & \cdots & {\beta}_{n-1} & {\beta}_n \\
                                                  \\
                                                  0 & {\beta}_0 & \cdots & {n-1 \choose 1}{\beta}_{n-2}& {n \choose 1}{\beta}_{n-1} \\
                                                  \\
                                                  0 & 0 & \cdots & {n-1 \choose 2}{\beta}_{n-3} & {n \choose 2}{\beta}_{n-2} \\
                                                  . & . & \cdots & . & . \\
                                                  . & . & \cdots & . & . \\
                                                  0 & 0 & \cdots & {\beta}_0 & {n \choose n-1}{\beta}_1
                                                  \end{array} \right|\end{array}$$\\

 $$\begin{array}{l}
 +\cdots+\frac{(-1)^{2n+1}x^{n-1}}{({\beta}_0)^{n+1}}\left| \begin{array}{ccccc}
                                                  {\beta}_0 & {\beta}_1 & {\beta}_2 & \cdots & {\beta}_n \\
                                                  \\
                                                  0 & {\beta}_0 & {2 \choose 1}{\beta}_1 & \cdots & {n \choose 1}{\beta}_{n-1} \\
                                                  \\
                                                  0 & 0 & {\beta}_0 & \cdots & {n \choose 2}{\beta}_{n-2} \\
                                                  . & . & . & \cdots & . \\
                                                  . & . & . & \cdots & . \\
                                                  0 & 0 & 0 & \cdots & {n \choose n-1}{\beta}_1
                                                  \end{array} \right| +\frac{x^n} {({\beta}_0)^{n+1}}\left| \begin{array}{ccccc}
                                                  {\beta}_0 & {\beta}_1 & {\beta}_2 & \cdots & {\beta}_{n-1} \\
                                                  \\
                                                  0 & {\beta}_0 & {2 \choose 1}{\beta}_1 & \cdots & {n-1 \choose 1}{\beta}_{n-2} \\
                                                  \\
                                                  0 & 0 & {\beta}_0 & \cdots & {n-1 \choose 2}{\beta}_{n-3} \\
                                                  . & . & . & \cdots & . \\
                                                  . & . & . & \cdots & . \\
                                                  0 & 0 & 0 & \cdots & {\beta}_0
                                                  \end{array} \right|. \end{array}\eqno(4.4)$$
\\

Since each minor in equation (4.4) is independent of $x$, therefore replacing $x$ by $\hat{M}_{LH}$ in equation (4.4) and then using the monomiality principle equation ${}_L{H}_n^{(m,r)}(x,y,z)=M_{LH}^n \{1\}~(n=1,2,3,\cdots)$, in the r.h.s. of the resultant equation, we find
$$ \begin{array}{l}
A_n(\hat{M}_{LH})=\frac{(-1)^n}{({\beta}_0)^{n+1}}\left| \begin{array}{ccccc}
                                                  {\beta}_1 & {\beta}_2 & \cdots & {\beta}_{n-1} & {\beta}_n \\
                                                  \\
                                                  {\beta}_0 & {2 \choose 1}{\beta}_1 & \cdots & {n-1 \choose 1}{\beta}_{n-2}& {n \choose 1}{\beta}_{n-1} \\
                                                  \\
                                                  0 & {\beta}_0 & \cdots & {n-1 \choose 2}{\beta}_{n-3} & {n \choose 2}{\beta}_{n-2} \\
                                                  . & . & \cdots & . & . \\
                                                  . & . & \cdots & . & . \\
                                                  0 & 0 & \cdots & {\beta}_0 & {n \choose n-1}{\beta}_1
                                                  \end{array} \right|\end{array}$$\\
\\
$$\begin{array}{l}
-\frac{(-1)^n{}_L{H}_1^{(m,r)}(x,y,z)}{({\beta}_0)^{n+1}}                                                  \left| \begin{array}{cccccc}
                                                  {\beta}_0 & {\beta}_2 & \cdots & {\beta}_{n-1} & {\beta}_n \\
                                                  \\
                                                  0 & {2 \choose 1}{\beta}_1 & \cdots & {n-1 \choose 1}{\beta}_{n-2}& {n \choose 1}{\beta}_{n-1} \\
                                                  \\
                                                  0 & {\beta}_0 & \cdots & {n-1 \choose 2}{\beta}_{n-3} & {n \choose 2}{\beta}_{n-2} \\
                                                  . & . & \cdots & . & . \\
                                                  . & . & \cdots & . & . \\
                                                  0 & 0 & \cdots & {\beta}_0 & {n \choose n-1}{\beta}_1
                                                \end{array} \right|\\

\\
+\frac{(-1)^n{}_L{H}_2^{(m,r)}(x,y,z) }{({\beta}_0)^{n+1}}\left| \begin{array}{cccccc}
                                                  {\beta}_0 & {\beta}_1 & \cdots & {\beta}_{n-1} & {\beta}_n \\
                                                  \\
                                                  0 & {\beta}_0 & \cdots & {n-1 \choose 1}{\beta}_{n-2}& {n \choose 1}{\beta}_{n-1} \\
                                                  \\
                                                  0 & 0 & \cdots & {n-1 \choose 2}{\beta}_{n-3} & {n \choose 2}{\beta}_{n-2} \\
                                                  . & . & \cdots & . & . \\
                                                  . & . & \cdots  & . & . \\
                                                  0 & 0 & \cdots & {\beta}_0 & {n \choose n-1}{\beta}_1
                                                  \end{array} \right|+\cdots+\end{array}$$\\

 $$\begin{array}{l}
 \frac{(-1)^{2n+1}{}_L{H}_{n-1}^{(m,r)}(x,y,z)}{({\beta}_0)^{n+1}}\left| \begin{array}{cccccc}
                                                  {\beta}_0 & {\beta}_1 & {\beta}_2 & \cdots & {\beta}_n \\
                                                  \\
                                                  0 & {\beta}_0 & {2 \choose 1}{\beta}_1 & \cdots & {n \choose 1}{\beta}_{n-1} \\
                                                  \\
                                                  0 & 0 & {\beta}_0 & \cdots & {n \choose 2}{\beta}_{n-2} \\
                                                  . & . & . & \cdots & . \\
                                                  . & . & . & \cdots & . \\
                                                  0 & 0 & 0 & \cdots & {n \choose n-1}{\beta}_1
                                                  \end{array} \right| +\frac{ {}_L{H}_n^{(m,r)}(x,y,z)} {({\beta}_0)^{n+1}}\left| \begin{array}{cccccc}
                                                  {\beta}_0 & {\beta}_1 & {\beta}_2 & \cdots & {\beta}_{n-1} \\
                                                  \\
                                                  0 & {\beta}_0 & {2 \choose 1}{\beta}_1 & \cdots & {n-1 \choose 1}{\beta}_{n-2} \\
                                                  \\
                                                  0 & 0 & {\beta}_0 & \cdots & {n-1 \choose 2}{\beta}_{n-3} \\
                                                  . & . & . & \cdots & . \\
                                                  . & . & . & \cdots & . \\
                                                  0 & 0 & 0 & \cdots & {n \choose n-1}{\beta}_1
                                                  \end{array} \right|. \end{array}\eqno(4.5)$$

Now, using operational rule (2.7) in the l.h.s. and combining the terms in the r.h.s. of equation (4.5), we obtain:
$$_{{}_LH^{(m,r)}}A_{n}(x,y,z)=\frac{{(-1)}^n}{{(\beta_0)}^{n+1}}\hspace{9cm}$$
$$\begin{array}{l}\left|\begin{array}{cccccc}
 1  & {}_L{H}_1^{(m,r)}(x,y,z) &  {}_L{H}_2^{(m,r)}(x,y,z)  & \cdots &  {}_L{H}_{n-1}^{(m,r)}(x,y,z)  & {}_L{H}_n^{(m,r)}(x,y,z) \\
 \\
 {\beta}_0  &  {\beta}_1  &  {\beta}_2  & \cdots &  {\beta}_{n-1}  &  {\beta}_n \\
 \\
 0  &  {\beta}_0  &  {2  \choose  1}{\beta}_1  & \cdots & {n-1 \choose 1}{\beta}_{n-2} & {n \choose 1}{\beta}_{n-1} \\
 \\
 0 & 0 & {\beta}_0 & \cdots & {n-1 \choose 2}{\beta}_{n-3} &  {n \choose 2}{\beta}_{n-2}\\
 . & . & . & \cdots & . & . \\
 . & . & . & \cdots & . & . \\
 0 & 0 & 0 & \cdots & {\beta}_0 & {n \choose n-1}{\beta}_1
 \end{array} \right|,\\
 {\beta}_n=-\frac{1}{A_0}\Big(\sum\limits_{k=1}^{n}{n \choose k}A_k{\beta}_{n-k}\Big),\\
 \hspace{4.35in}n=1,2,3,\cdots,\end{array}\eqno(4.6)$$
where ${\beta}_0,~{\beta}_1,\cdots,{\beta}_n \in \mathbb{R},~{\beta}_0 \neq 0$ and ${}_L{H}_n^{(m,r)}(x,y,z)~(n=0,1,\cdots)$ are the Laguerre-Gould Hopper polynomials defined by equation (1.1).\\

Equation (4.3) together with equation (4.6), gives the determinantal definition for the Laguerre-Gould Hopper based Appell polynomials (LGHAP) $_{{}_LH^{(m,r)}}A_{n}(x,y,z)$.\\

Also, by giving suitable values to the variables and indices we can find the determinantal definitions for the members belonging to the Laguerre-Gould Hopper based Appell family.

\end{document}